
\input gtmacros
\input rlepsf
\input newinsert
\input gtmonout
\volumenumber{1}
\volumeyear{1998}
\volumename{The Epstein birthday schrift}
\pagenumbers{479}{492}
\received{20 November 1997}
\revised{7 November 1998}
\published{17 November 1998}
\papernumber{23}

\def\ORR{{\overline R}}
\def\nek{,\ldots,}
\def\calF{{\cal F}}
\def\calG{{\cal G}}

\def\calP{{\cal P}}
\def\Def{\mathop{\rm Def}}
\def\s3{{\bf S^3}}

\let\bet\beta
\let\Gam\Gamma
\let\gam\gamma
\let\eps\varepsilon
\let\sig\sigma
\let\Sig\Sigma
\let\tet\theta
\let\Lam\Lambda
\let\Ome\Omega
\let\bks\backslash
\let\part\partial

\let\sect\section
\let\subheading\rk
\let\pr\prf
\let\qua\stdspace

\reflist

\refkey\A
{\bf L\,V~Ahlfors}, {\it Finitely generated Kleinian groups},
Amer. J. Math. 86 (1964) 413--429; 87 (1965) 759

\refkey\Be
{\bf L~Bers}, {\it On boundaries of Teichm\"uller spaces and on Kleinian
groups I}, Annals of Math. 91 (1970) 570--600

\refkey\Bo
{\bf B~Bowditch}, {\it Geometrical finiteness of hyperbolic groups},
preprint, University of Melbourne

\refkey\BoM
{\bf B~Bowditch}, {\bf G~Mess}, {\it A 4--Dimensional Kleinian Group},
Trans. Amer. Math. Soc. 344 (1994) 391--405

\refkey\J
{\bf T~J\o rgensen}, {\it Compact $3$--manifolds of constant
negative curvature fibering over the circle}, Annals of Math. 
106 (1977) 61--72

\refkey\K
{\bf M~Kapovich}, {\it On Absence of Sullivan's cusp finiteness
theorem in higher dimensions}, preprint

\refkey\KP {\bf M~Kapovich}, {\bf L~Potyagailo}, {\it On absence of
Ahlfors' finiteness theorem for Kleinian groups in dimension $3$},
Topology and its Applications, 40 (1991) 83--91

\refkey\KPS {\bf M~Kapovich}, {\bf L~Potyagailo}, {\it On absence of Ahlfors'
and Sullivan's  finiteness theorems for Kleinian groups in higher
dimensions}, Siberian Math. Journal 32 (1992) 61--73

\refkey\L {\bf D~Long}, {\it Immersions and embeddings of totally geodesic
surfaces}, Bull. London Math. Soc. 19 (1987) 481--484

\refkey\Ma
{\bf B~Maskit}, {\it On boundaries of Teichm\"uller spaces and on
Kleinian groups, II}, Annals of Math. 91 (1970) 608--638

\refkey\Mab
{\bf B~Maskit}, {\it On Klein's Combination theorem III},
from: ``Advances in the theory of Riemann Surfaces'', Princeton
Univ. Press (1971) 297--310

\refkey\Mac
{\bf B~Maskit}, {\it Kleinian groups}, Springer--Verlag (1988)

\refkey\Mo
{\bf J~Morgan}, {\it Group action on trees and the compactification
of the space of conjugacy classes of $SO(n,1)$--representations},
Topology 25 (1986) 1--33

\refkey\P
{\bf L~ Potyagailo}, {\it Finitely generated Kleinian groups in 3--space
and 3--manifolds of infinite
homotopy type}, Trans. Amer. Math. Soc. 344
(1994)  57--77

\refkey\Po
 {\bf L~Potyagailo}, {\it The problem of finiteness for Kleinian
groups in $3$--space}, from: ``Proceedings of International Conference,
Knots-90", Osaka (1992)

\refkey\Sc {\bf P~Scott}, {\it Finitely generated $3$--manifold groups
are finitely presented}, J. London Math. Soc. 6  (1973)
437--440

\refkey\Scb {\bf P~Scott}, {\it Subgroups of surface groups are almost
geometric}, J. London Math. Soc. 17 (1978) 555--565; Correction
ibid 32  (1985) 217--220

\refkey\Su
{\bf D~Sullivan}, {\it Travaux de Thurston sur les groupes
quasi-fuchsiens et les varietes hyperboliques de dimension $3$
fibres sur $S^1$}, Lecture Notes in Math, 842, Springer--Verlag,
Berlin--New York (1981) 196--214

\refkey\T
{\bf W~Thurston}, {\it The geometry and topology of
$3$--manifolds}, Princeton University Lecture Notes  
(1978)

\endreflist

\title{The boundary of the deformation space of the\\\vglue-2.2mm\\fundamental
group of some hyperbolic\\3--manifolds fibering over the circle}
\covertitle{The boundary of the deformation space of the\\fundamental
group of some hyperbolic\\3--manifolds fibering over the circle}
\asciititle{The boundary of the deformation space of the\\fundamental
group of some hyperbolic\\3-manifolds fibering over the circle}
\shorttitle{Deformation space of hyperbolic 3--manifolds}

\author{Leonid Potyagailo}

\address{D\'epartement de Math\'ematiques\\Universit\'e de Lille
1\\
59655 Villeneuve d'Ascq, France}

\asciiaddress{Departement de Mathematiques\\Universite de Lille 1\\
59655 Villeneuve d'Ascq, France}

\email{potyag@gat.univ-lille1.fr}

\abstract
By using Thurston's bending construction  we obtain  a
sequence of faithful discrete representations $\rho_n$ of the
fundamental group of a closed hyperbolic $3$--manifold fibering over
the circle into the isometry group $Iso\ {\bf H}^4$ of the
hyperbolic space ${\bf H}^4 $. The algebraic limit of $\rho_n$
contains a finitely generated subgroup $F$ whose $3$--dimensional
quotient $\Omega(F)/F$ has infinitely generated fundamental group,
where $\Omega(F)$ is the discontinuity domain of $F$ acting on the
sphere at infinity $S^3_{\infty}=\partial{\bf H}^4 $. Moreover $F$
is isomorphic to the fundamental group of a closed surface and
contains infinitely many conjugacy classes of maximal parabolic
subgroups. 
\endabstract

\asciiabstract{%
By using Thurston's bending construction  we obtain  a
sequence of faithful discrete representations \rho_n of the
fundamental group of a closed hyperbolic 3-manifold fibering over
the circle into the isometry group Iso H^4 of the
hyperbolic space H^4. The algebraic limit of \rho_n
contains a finitely generated subgroup F whose 3-dimensional
quotient \Omega(F)/F has infinitely generated fundamental group,
where \Omega(F) is the discontinuity domain of F acting on the
sphere at infinity. Moreover F
is isomorphic to the fundamental group of a closed surface and
contains infinitely many conjugacy classes of maximal parabolic
subgroups.}


\maketitle

\sect{Introduction and statement of results}

By a Kleinian (discontinuous) group $G$ we mean a subgroup of the
group ${\rm Conf}({\bf S}^n)\cong SO_+(1,n+1)$  of conformal
transformations of $\ORR^{\,n}=S^n=R^n\cup\{\infty\}$ which acts
discontinuously on a non-empty set $\Ome (G)\subset S^n$ called its
domain of discontinuity. It may be connected or not; we will say
that $G$  is a function group if there is a connected component
$\Ome_G\subset \Ome (G)$ that is invariant under the action of the
whole group: $G\Ome_G=\Ome_G$.  The quotient spaces $M_G=\Ome_G/G$
and $M(G)=\Ome(G)/G$  are $n$--manifolds in  the case  in which $G$
is torsion-free.  The complement $\Lam (G)=(S^n\bks\Ome
(G))\subset \part {\bf H}^{n+1}$  is called the limit set of $G$.

A finitely generated Kleinian group $G$  is called geometrically
finite if for some $\eps
>0$ there exists an $\eps$--neighbourhood of $H_G/G$ in ${\bf
H}^{n+1}/G$ which is of finite hyperbolic volume. Here
$H_G\subset{\bf H}^{n+1}$  is the convex hull of $\Lam (G)$.

Let us consider for $n=3$ a hyperbolic $3$--manifold $M=H^3/\Gam$
 ($\Gam\subset
PSL_2{\bf C}$)  fibering over the circle $S^1$  with fiber a
closed surface $\sig$. The notation is $\ M=\sig\tilde\times S^1$.
A representation $\rho\co\pi_1(M)\to {\rm Conf}({\bf S}^3)$   is
called admissible if the following conditions are satisfied.
\items
\item{(1)} $\rho\co\Gam\to {\rm Conf}({\bf S}^3)$  is faithful and $\rho(\Gam)=\Gamma_0$
is Kleinian.
\item{(2)} $\rho$  preserves the type of each element, ie $\rho
(\gam)$  is loxodromic for all $\gam\in \Gam$.
\item{(3)} $\rho$  is induced by a homeomorphism
$f_\rho\co\Ome(\Gam )\to \Ome (\Gam_0)$, namely $f_\rho \gam
f_\rho^{-1}=\rho (\gam)$, $\gam \in \Gam$.\enditems

The set of all admissible representations modulo conjugation in
${\rm Conf}({\bf S}^3)$  is called the deformation space
Def$(\Gam)$  of the group $\Gam$.

The set Def$(\Gam)$ inherits the topology of convergence on
generators of $\Gam$  on compact subsets in ${\bf S}^3$ because
 Def$(\Gam )\subset \left({\rm Conf}({\bf S}^3)\right)^k/\!\sim$, $k\in
{\bf N}$ ($\sim$ is  conjugation in ${\rm Conf}({\bf S}^3)$). As
Def$(\Gam)$ is a bounded domain [\Mo] two questions have arisen.
The first is to describe the cases when Def$(\Gam)$ is non-trivial
and the second is to study the boundary $\part\Def(\Gam)$,  as
 was done for the classical Teichm\"uller space [\Be], [\Ma].
The answer to the first question is still unknown even in the
case when $M$ is Haken.  We will consider the case when $M$
contains many totally geodesic surfaces. Each of them
produces a curve in Def$(\Gam)$ by Thurston's ``bending"
construction [\T].
Our main interest is in groups which appear on the
boundary $\part\Def (\Gam)$.  These are higher dimensional
analogs
of $B$--groups which arise as the limits of sequences of
quasifuchsian groups in classical Teichm\"uller space.

One of the most fundamental questions  is to describe the
topological type of the orbifold $M(\Gamma)=\Omega (\Gam)/\Gamma$
(a manifold in the case when $\Gamma$ is torsion-free), in
particular, when $\Gamma$ is a function group it is important to
know when the fundamental group $\pi_1 (M_G =\Omega_\Gam /\Gam)$
turns out to be finitely generated, or even more generally when it
has finite homotopy type.

In dimension $2$ the famous  theorem of Ahlfors [\A] says that a
finitely generated non-elementary Kleinian group $G\subset {\rm
Conf}({\bf R}^2)$ has a factor-space $\Ome (G)/G$ consisting of a
finite number of Riemann surfaces $S_1\nek S_n$ each having a
finite hyperbolic area.

We discovered in [\KP] that the weakest topological version of
Ahlfors' theorem  does not hold starting already with dimension
$3$. Namely we constructed a finitely generated function group
$F\subset {\rm Conf}({\bf S}^3)$ such that the group $\pi_1
(\Omega_F /F)$ is not finitely generated. Afterwards it was
pointed out in [\Po] that this group is in fact not finitely
presented.

It has also been shown that there exists a finitely generated
Kleinian group with infinitely many conjugacy classes of
parabolics [\K].

In [\P] we constructed a finitely generated group $F_1$ such that $\pi_1
(\Omega_{F_1} /F_1 )$ is not finitely generated and having
infinitely many non-conjugate elliptic elements; moreover $F_1$
appears as an infinitely presented subgroup of a geometrically
finite Kleinian group in ${\bf H}^4$ without parabolic elements.
On the other hand, it was shown in [\BoM] that a finitely
generated but infinitely presented group can also appear as a
subgroup of a cocompact group in $SO(1,4)$.

\proclaim {Theorem 1}  Let  $\Gam =\pi_1(M)$  be the fundamental
group of a hyperbolic $3$--manifold $M$ fibering over the circle
with fiber a closed surface $\sigma$. Suppose that $\Gamma$ is
commensurable with the reflection group $R$ determined by the
faces of a right-angular polyhedron $D\subset{\bf H}^3$.  Then
there exists a finite-index subgroup $L\subset\Gamma$ and a path
 $\beta_t\co[0,1[\mapsto \Def(\Gamma)$ such that $\beta_t$
 converges to a faithful representation
 $\beta_1\in\partial{\rm Def}(\Gamma)$ (as $t\to 1$)
 and the following hold:

\items
\item{\rm (1)} $\beta_1(F_L)$ contains
infinitely many conjugacy classes of maximal parabolic subgroups,
\item{\rm (2)} $\displaystyle\pi_1(\Omega_{\beta_1(F_L)})/\beta_1
(F_L)$ is  infinitely generated,
\enditems

where $F_L=L\cap\pi_1\sigma$ is isomorphic to the fundamental group of
a closed hyperbolic surface which finitely covers $\sigma$ and
$\beta_1(F_L)$ acts discontinuously on an invariant component
$\Omega_{\beta_1(F_L)}\subset {\bf S}^3$.
\endproc

 \subheading{Remark}  Groups satisfying all the conditions of
Theorem 1 do exist.  An example of Thurston, of the reflection
group in the faces of the right-angular dodecahedron, which is
commensurable with a group of a closed surface bundle, is given in
[\Su].

\subheading{Acknowledgement} This paper was prepared several
years ago while the author had a Humboldt Fellowship at the
R\"uhr-Universit\"at in Bochum. The author is deeply grateful to
Heiner Zieschang and to the Humboldt Foundation for this
opportunity. I would also like to thank Nicolaas Kuiper (who died
recently) for reading a preliminary version of the manuscript and to
express my gratitude to the referee for many useful remarks and
corrections.

\sect{Outline of the proof}

Before giving a formal proof of the Theorem let us describe it
informally.

Our construction is inspired essentially by  papers [\K], [\KPS]
and [\P].  In the first two a free Kleinian group of finite rank
satisfying the conclusion $(2)$ was produced, whereas now we give
an example of a closed surface group with this property. Our
present construction is essentially easier than that of [\P]. Also,
we produce a curve in the deformation space whose limit point is
the group in question.

\rk{Step 1}We start with  an uniform lattice $\Gam\subset
PSL_2{\bf C}$ commensurable with the reflection group $R$ whose limit
set is the Euclidean $2$--sphere $\partial B_1$ -- the boundary of the
ball $B_1\subset \s3$. There exists a Fuchsian subgroup
$H_2\subset\Gam$ leaving invariant a vertical plane $\pi$ whose
intersection with $B_1$ is a round circle, its limit set $\Lambda
(H_2)$ (see figure 1). The group $H_2$ also leaves invariant a geodesic
plane $w_2\subset B_1$.  Consider the action of the group $\Gam$ in
the outside ball $B_1^*=\s3\setminus B_1$. For some finite-index
subgroup $\Gam_1$ of $\Gam$ we construct a new group $G_1$ obtained by
Maskit's Combination theorem from $\Gamma_1$ and
$\tau_{\pi}\Gamma_1\tau_{\pi}$ combined along the common subgroup
$H_2={\rm Stab\ }w_2$, where $\tau_{\pi}$ is the reflection in
$\pi$. The new group $G_1$ is still isomorphic to some subgroup
$G^*\subset R$ of finite index essentially because the same
construction can be done inside $B_1$ by reflecting the picture
along the geodesic plane $w_2$. Thus $G_1$ belongs to the deformation
space ${\rm Def}(G^*_1)$. One can obtain a fundamental domain
$R(G_1)\subset B_1^*$ of $G_1$ which is situated in a small
neighbourhood of the spheres $\partial B_1$ and $\tau_{\pi}(\partial
B_1)$.

\fig{1}\relabelbox\small
\epsfxsize 3.5truein\epsfbox{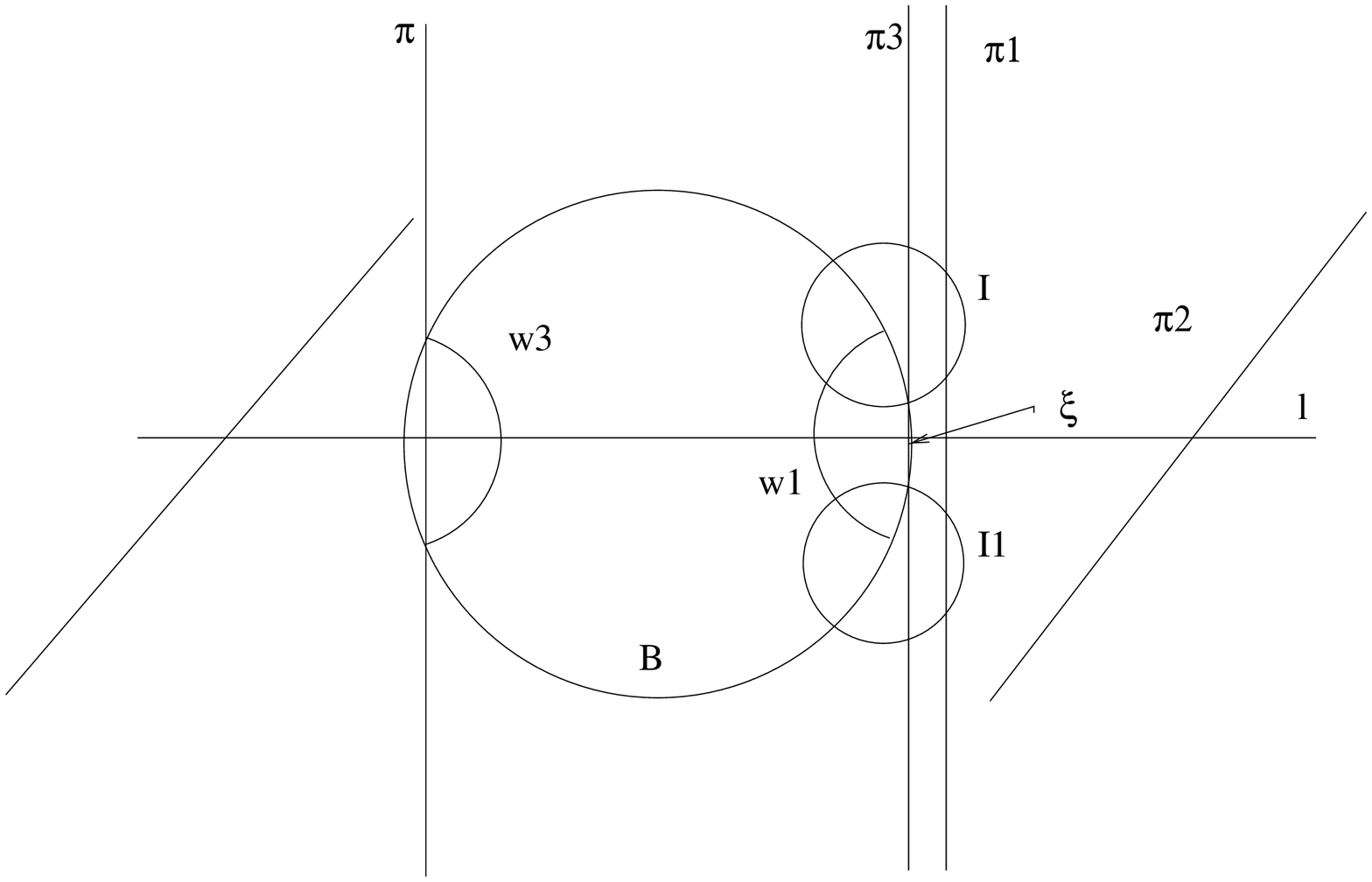}
\adjustrelabel  <-2pt,-2pt> {w3}{$w_2$}
\adjustrelabel  <-2pt,0pt> {w1}{$w_1$}
\relabel {B}{$B_1$}
\relabel {I}{$I_{g_1}$}
\relabel {I1}{$I'_{g_1}$}
\relabel {l}{$\ell$}
\adjustrelabel  <-4pt,-2pt> {p}{$\pi$}
\adjustrelabel  <-5pt,-0pt> {p3}{$\pi_3$}
\adjustrelabel  <-4pt,2pt> {p1}{$\pi_1$}
\adjustrelabel  <2pt,-2pt> {p2}{$\pi_2$}
\relabel {x}{$\xi$}
\endrelabelbox
\endfig 

\fig{2}\relabelbox\small
\epsfxsize 3.5truein\epsfbox{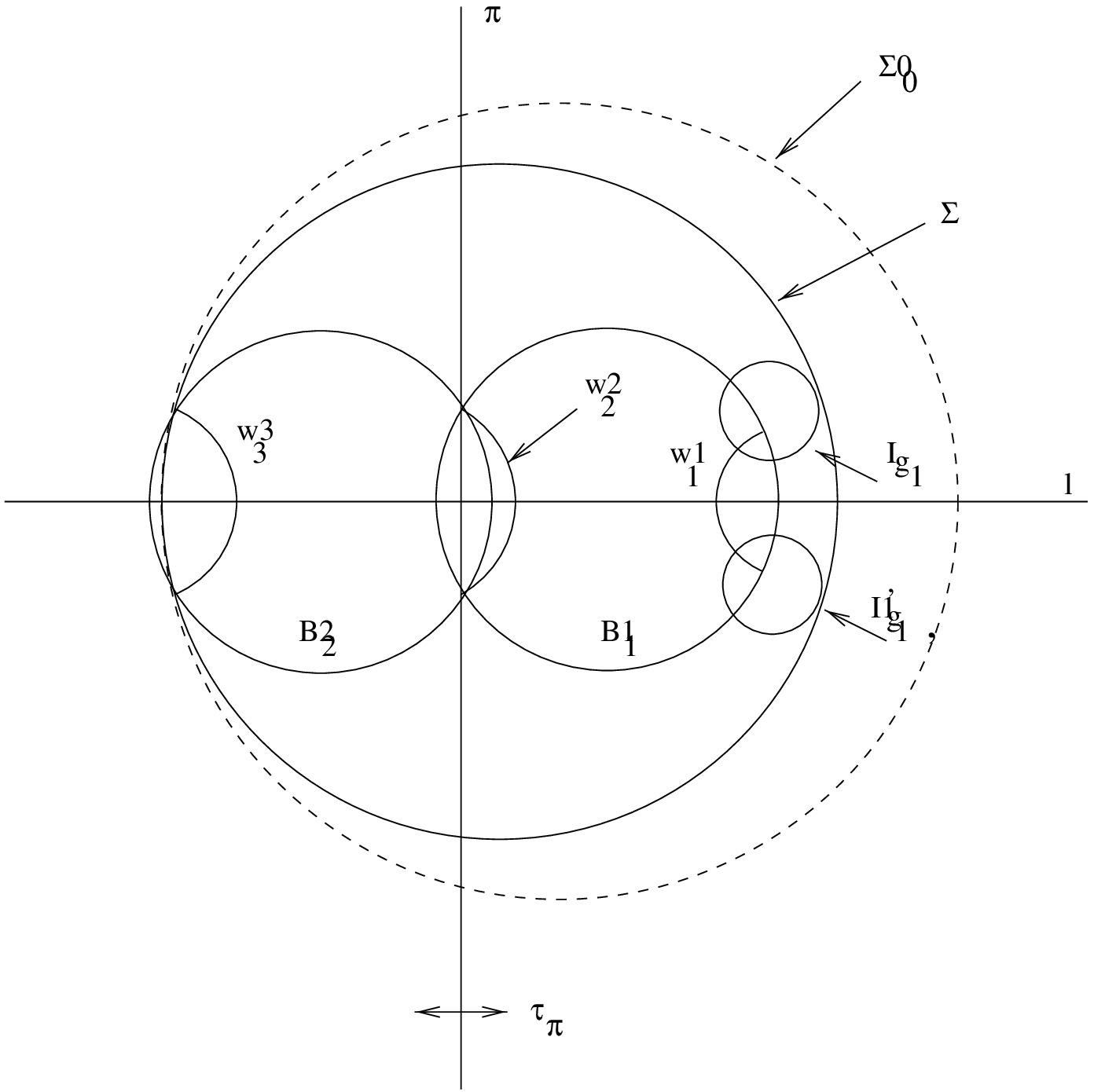}
\adjustrelabel <-2pt, -2pt> {w3}{$w_3$}
\adjustrelabel <-2pt, -2pt> {w2}{$w_2$}
\relabel {w1}{$w_1$}
\relabel {B1}{$B_1$}
\relabel {B2}{$B_2$}
\adjustrelabel <-1pt, -4pt> {I}{$I_{g_1}$}
\adjustrelabel <0pt, -3pt> {I1}{$I'_{g_1}$}
\relabel {l}{$\ell$}
\relabel {p}{$\pi$}
\relabel {t}{$\tau_{\pi}$}
\adjustrelabel <-1pt, 0pt> {S}{$\Sigma$}
\adjustrelabel <-2pt, 0pt> {S0}{$\Sigma_0$}
\endrelabelbox
\endfig 

\rk{Step 2}There is another geodesic plane $w_1\subset B_1$ disjoint
from $w_2$ whose stabilizer in $\Gamma_1$ is $H_1$ (see figure 2).
Denote by $B_2$ the ball $\tau_\pi(B_1)$.
Take a sphere $\Sigma\subset B_1^*$ passing through
the circle $w_3\cap B_2$ -- the limit set of the group $\tau_{\pi}
H_1\tau_{\pi}$ -- and tangent to the isometric spheres of some
element $g_1\in\Gam_1$, where $H_1$ is a subgroup of $\Gam_1$
stabilizing $w_1$. We now construct a family of Euclidean spheres
$\Sigma_t\ (0\leq t\leq 1,\ \Sigma_1=\Sigma)$
and corresponding groups $\calG_t$ obtained as
before from $G_1$ and
$\displaystyle\tau_{\Sigma_t}G_1\tau_{\Sigma_t}$ by using the
combination method along common closed surface subgroups. We prove
then that there is a path $\beta_t \co t\in [0,1[\ \mapsto
\beta\in\Def (L')$ such that $\beta_0=L',\ \beta_t=\calG_t$ where
$L'$ is some finite-index subgroup of $R$. One can equally say
that $\beta_t$ is obtained by using Thurston's bending
deformation. The main point is now to prove that the limit group
$\displaystyle \calG_1=\lim_{t\to 1}\beta_t(L')$ is discontinuous
and has a fundamental domain obtained from the part of $R(G_1)$ by
doubling along the sphere $\Sigma$. The group $\calG_1$ is also
isomorphic to $L'$  and so contains a fundamental group $\cal N$
of a closed surface bundle over the circle which is isomorphic to
the group $L=\Gamma\cap L'$. Let $\cal F$ be the fundamental group
of the fiber  given by $\beta_1(F_L=F\cap L)$. Since two isometric
spheres of the element $g_1\in\Gamma_1$ are tangent to $\Sigma$,
we get a new accidental parabolic element $g=g_1\cdot g_2,\
g_2=\tau_\Sigma g_1\tau_\Sigma$ in the group $\calG_1$. By a
choice of $g_1$ made from the very beginning we assure that
$g\in\cal F$, so we have a pseudo-Anosov action of some element
$t\in\cal N\setminus\cal F$ such that the orbit $t^n\cdot g\cdot
t^{-n}\ (n\in {\bf Z})$ gives us infinitely many conjugacy classes
of maximal parabolic subgroups of $\cal F$. Now Scott's compact
core theorem implies that $\pi_1(\Ome_{\cal F})/{\cal F}$ is not
finitely generated. \hfill {\it End of outline}

\sect{Preliminaries}

 We will consider the Poincar{\'e} model of
hyperbolic space ${\bf H}^3$ in the unit ball $B_1$  equipped with
the hyperbolic metric $\rho$. By a right-anguled polyhedron
$D\subset {\bf H}^3$  we mean a polyhedron all of whose dihedral
angles are $\pi/2$.

Consider the tesselation of ${\bf H^3}$ by images of $D$ under
the reflection group
$R$ from Theorem 1. Denote by $W\subset{\bf H}^3$ the collection
of geodesic planes $w$ such
that    there exists $r\in R$,  for which
$r(w)\cap\part D$  is a face of $D$.

It is easy to see that if $\sig_1$ and $\sig_2$ are two faces of
$D$ with
 $\sig_1\cap \sig_2=\emptyset$,  then also the  geodesic planes
$\tilde\sig_1\supset\sig_1$ and $\tilde\sig_2\supset\sig_2$ have
no point in common.  One can easily show that the distance between
$\sig_1$ and $\sig_2$, as well as that of $\tilde\sig_1$ and
$\tilde\sig_2$, is realized by a common perpendicular $\ell$ for
which $\ell\cap {\rm int} D\not=\emptyset$.

Let $\Gam_0=R\cap\Gam$  which is a subgroup of a finite index in
both groups $R$  and $\Gam$.  By passing to a subgroup of a finite
index and preserving notation, we may assume that $\Gam_0$  is a
normal subgroup in $R$, $|R:\Gam_0|<\infty$. For a plane $w\in W$  we
write $H_w={\rm Stab}(w,\Gam_0)=\{ g\in \Gam_0,\ gw=w\}$.  It is
not hard to see that $H_w$  is a Fuchsian group of the first kind
commensurable with the reflection group determined by the edges of
some face of the polyhedron $r(D_1),\ r\in R$.

Let us now fix two disjoint planes $w_1$ and $w_2$ from $W$
containing opposite faces of $D$ and let $\ell$ be their common
perpendicular; up to conjugation in ${\rm Isom\ } {\bf H}^3$ we
can assume that $\ell$ is a Euclidean diameter of $B_1$. Denote
$B^*_1={\bf S}^3\bks cl(B_1)$ as well (where $cl (\cdot )$  is the
closure of a set). We have the following:

\proclaim{Lemma 1} For every horosphere $\pi_3$ in $B_1^*$
centered at the point $\xi\in \ell\cap\partial B_1$ (see  figure 1)
there exists $\eps_0>0$ such that for every $\eps$--close sphere
$\pi_1\subset B_1^*$  to $\pi_3$ ($\eps<\eps_0$) orthogonal to the
plane $\pi_2$ there exists a
geodesic plane $w$ and an element $g_1\in [H_w, H_w]$ (commutator
subgroup) such that:
$$\eqalign{I_{g_1}\cap\pi_1\not=\emptyset \ \  {\rm and} \ \
g_1(I_{g_1}\cap\pi_1)=I'_{g_1}\cap\pi_1,\cr {\rm where}\quad
I_{g_1},\ I'_{g_1}=I_{g_1^{-1}} \quad\hbox{are isometric spheres
of $g_1$}.}\eqno(1)$$\endproc

\pr Up to further conjugation in ${\rm Isom~}B_1$ preserving $\ell$
we may assume that $\pi_3$ is the vertical plane tangent to
$\partial B_1$ at $\xi\in \ell\cap\partial B_1$. Take $w=w_1$ and let
$g_1\in [H_{w_1}, H_{w_1}]$ be any primitive element corresponding
to a simple dividing loop on the surface $w_1/H_{w_1}$.

Suppose first that $I_{g_1}\cap \pi_3=\emptyset$. In this case we
proceed as follows. Put $\chi =\tau_{w_1}\circ\tau_{w_2}\in R$,
where $\tau_{w_i}$  denotes the reflection in plane $w_i\ (i=1,
2)$. Then $\chi$ is a hyperbolic element whose invariant axis is
$\ell$. Consider  the sequence of planes $\chi^n(w_1)$.  We claim that,
for some $n$, $\chi^n(I_{g_1})\cap\pi_3\not=\emptyset$. In fact
this follows directly from the fact that  the fixed point $\xi$ of
the hyperbolic element $\chi$ is a conical limit point of
$\Gam_0$, and so the approximating sequence $\chi^n(I_{g_1})$
should intersect a fixed horosphere (or equivalently by sending
$\xi$ to the infinity and passing to the half-space model one can
see that $\chi$ becomes now a dilation $z\mapsto \lambda z \
(\lambda
> 0)$ which implies that the translations of the image of
$I_{g_1}$ by powers of the dilation will intersect a fixed
horosphere at infinity). Since $\Gam_0$ is normal in $R$ it now
follows that $\chi^n g_1\chi^{-n}\in [H_{\chi^n(w_1)},
H_{\chi^n(w_1)}]\subset\Gam_0$ and $\displaystyle\chi^n(I_{g_1}) =
I_{\chi^n g_1\chi^{-n}}$. The latter is true since $\chi$
preserves each Euclidean plane passing through $B_1\cap \ell$ and,
hence $(\chi^n g_1\chi^{-n})\vert_{\chi^n(I_{g_1})}$is an Euclidean
isometry. So up to replacing $w_1$ by $\chi^n(w_1)$ and $g_1$ by
$\chi^n g_1\chi^{-n}$ if needed, we may assume that
$I_{g_1}\cap\pi_3\not=\emptyset$. The same conclusion is then
obviously  true for a plane $\pi_1\subset B_1^*$ sufficiently
close to $\pi_3$.

For $\ell_1=I_{g_1}\cap\pi_1$ we now claim that
$g_1(\ell_1)=\ell_2=I'_{g_1}\cap\pi_1$.  Indeed,
$g_1=\tau_{\pi_2}\cdot\tau_{I_{g_1}}$  where $\pi_2$  is
orthogonal to $\pi_1$  and contains $\ell$  (figure 1). Evidently
$$g_1(\ell_1)=\tau_{\pi_2}\left(I_{g_1}\cap\pi_1\right)
=\tau_{\pi_2}(I_{g_1})\cap\pi_1=I'_{g_1}\cap\pi_1\hfill\eqno(2)$$
since $\tau_{\pi_2}(\pi_1)=\pi_1$. The lemma is proved.
\endproof

So we can suppose  that $w_1\in W$ is chosen satisfying all the
conclusions of Lemma 1. Let $w_2\in W$ be a geodesic plane
disjoint from $w_1$ and let
 $\ell$ be their common
perpendicular  passing through the origin of $B_1$. Now consider
the Euclidean plane $\pi$ orthogonal to $\ell$ (figure 2) such
that $$\pi\cap\part B_1=\pi\cap w_2\ .$$ It is not hard to see
that ${\rm Stab}(\pi, \Gam)= {\rm Stab}(w_2,\Gam)=H_{w_2}$. Reflecting our
picture in the plane $\pi$ we get $$\eqalign{B_2&=\tau_\pi (B_1)\
,\quad w_3=\tau_\pi (w_2)\quad {\rm and}\cr H_{w_3}&=\tau_\pi
H_{w_1}\tau_\pi\ .\cr}$$

By Lemma 1 we can now find  a Euclidean sphere $\Sig$ centered on
$\ell$ which goes through the circle $w_3\cap\part B_2$ and is
tangent to $I_{g_1}$ (figure 2).  Moreover, by Lemma 1, $\Sig$ is
tangent also to $I'_{g_1}$.

Denote $\Sig'=\tau_\pi^{-1}(\Sig)$.

\proclaim{Lemma 2}  There exists a subgroup $\Gam_1\subset\Gam_0$
of  finite index such that the following conditions hold:

\items
\item{\rm (a)} The boundary of the isometric fundamental domain
$\calP (\Gam_1)\subset B^*_1$  lies in a regular
$\eps$--neighbourhood of $\part B_1^* \left(B^*_1={\bf S}^3\bks cl
(B_1),\ \eps
>0\right)$.

\item{\rm (b)} $\Sig\cap I_\gam =\emptyset\ ,\quad \gam \in
\Gam_1\bks
\{g_1,g_1^{-1}\}$.

\item{\rm (c)} For subgroups $H_1=\Gam_1\cap H_{w_1},
H_2=\Gam_1\cap H_{w_2}$  there exists another fundamental domain
$R(\Gam_1)\subset B^*_1$  of  $\Gam_1$  such that $$R(\Gam_1)\cap
(\pi\cup\Sig')=\calP (H)\cap (\pi \cup\Sig'),$$ where $\calP(H)$
is an isometric fundamental domain for the group $H=\langle
H_1,H_2\rangle$.

\item{\rm (d)} $g_1\in \Gam_1\cap [H_1,H_1]$.\enditems\endproc

\pr  This Lemma can be obtained by repeating the arguments of [\P,
Main Lemma]. We just sketch these considerations.  First, we
choose
a subgroup $\tilde\Gamma \subset\Gamma_0$  of a finite index
satisfying conditions (a) and (b) such that $g_1\in\tilde\Gamma$
 by
using the property of separability of infinite cyclic subgroups in
$\Gam_0\ [\L]$.

To obtain (c) we will find $\Gam_1$ by using Scott's
$LERF$--property of the group $\Gam_0$ with respect to its
geometrically finite subgroups (see [\Sc], [\Scb]). To this end we
proceed as follows: the group $H$ is geometrically finite as a
result of Klein--Maskit free combination from $H_1$ and $H_2$, which
are both geometrically finite subgroups of $\Gam_0$. The $LERF$
property now says that  for the element $g_1$ there exists a
subgroup of $\Gam_0$ of finite index which contains $H$ and does
not contain $g_1$. Call this subgroup $\Gam_1$. Evidently, $g_1\in
[H_1,H_1]\subset\Gam_1$ by construction. For the complete proof,
see [\P, Main Lemma].
\endproof

Let us introduce the following notation:
$\Ome^-_1=B^*_1\bks\bigcup_{\gam\in \Gam_1}\gam (\pi^-)$  where
$\pi^-$  is the component of ${\bf S}^3\bks\pi$  for which
$w_3\in\pi^-$.  Let $\Gam'_1={\rm Stab}(\Ome^-_1,\Gam_1)$.

The complete proof of the following assertion can be also
found
in [\P, Lemma~3].

\proclaim{Lemma 3}  The group $G_1=\langle \Gam'_1,\
\tau_\pi\Gam'_1\tau_\pi\rangle$  is discontinuous and

\items
\item{\rm(1)} $G_1\buildrel\sim\over =\Gam'_1*_{H_2}
\left(\tau_\pi\Gam'_1\tau_\pi\right)$.

\item{\rm(2)}  $G_1$  is isomorphic to a subgroup
$G^*_1\subset R$  of  finite index.\enditems\endproc

\proof{Sketch of proof} {(1)}\qua This
follows from the fact that the plane $\pi$
is strongly invariant under $H_2$ in $\Gam'_1$  by [\P, Lemma
3.c], which means $H_2\pi=\pi$  and
$\gam\pi\cap\pi=\emptyset\ ,\ \gam\in\Gam'_1\bks H_2$.
One can now get assertion (1) from Maskit's First Combination theorem
[\Mab].

(2)\qua Consider the reflection $\tau_{w_2}$ in the
geodesic plane $w_2\subset B_1$.  We claim that the group
$G^*_1=\langle
\Gam'_1,\tau_{w_2}\Gam'_1\tau_{w_2}\rangle$  is isomorphic to
$G_1$.  Indeed, $w_2$  is also strongly invariant under $H_2$
in $\Gam'_1$  and we again observe that
$G^*_1=\Gam'_1*_{H_2}\left(\tau_{w_2}\Gam'_1\tau_{w_2}\right)
\buildrel\sim\over = G_1$  because
$\tau_{w_2}\mid_{{}_{w_2}}=\tau_\pi \mid_{{}_{_\pi}}
=id$.

Now $\tau_{w_2}\in R$.  Therefore, $G^*_1\subset R$ and $G^*_1$
has a compact fundamental domain $R(G^*_1)=R(\Gam'_1)\cap
\tau_{w_2}(R(\Gam'_1))$.  The covering ${\bf
H}^3\big/(G^*_1\cap\Gam_0)\to {\bf H}^3\big/G^*_1$ is finite since
$|R:\Gam_0|<\infty$  and, hence, the manifold $M\left(G^*_1\cap
\Gam_0\right)={\bf H}^3\big/(G^*_1\cap\Gam_0)$  is compact.  Thus,
the covering $M(G^*_1\cap\Gam_0)\to M(\Gam_0)$  is finite as well
and so $|\Gam_0:G^*_1\cap\Gam_0\mid <\infty$.
\endproof

\proclaim{Corollary 4}  There exists a path $\alpha_t \co[0,1]\ \to
Def(G^*_1)$ such that $\alpha_0=G^*_1$  and
$\alpha_1=G_1$.\endproc

\pr By choosing a continuous family of spheres $\mu_t$  for which
$\mu_t\cap\pi=w_2\cap\pi=\Lambda(H_2),\ \mu_0\supset w_2,\
\mu_1=\pi,\ t\in [0,1)$,  we construct the family of groups
$G_t=\langle\Gam'_1,\tau_{\mu_t}\Gam'_1\tau_{\mu_t} \rangle$  by
the arguments of Lemma 3.  Consider now the action of $\Gam'_1$ in
$B^*_1$ where $\ p_1\co B^*_1\to B^*_1/\Gam_1$  is the covering map.
The surfaces $p_1(\mu_t)$  are all embedded and parallel due to
condition (b).  If now $\Omega_{G_t}$  is the component of $G_1$
containing $\infty$  then the manifold $M_{G_t}=\Ome_{G_t}/G_t$ is
homeomorphic to the double of the manifold $M^-_1=
\Ome^-_1/\Gam'_1$ along the boundary $p_1(\pi)$.  Thus, for all
$t\in [0,1]$,  $M_{G_t}$  are all homeomorphic and there exists a
continuous family of homeomorphisms $f_t\co\Ome (G^*_1)\to\Ome
(G_t)$  such that $G_t=f_tG^*_1f_t^{-1}$,
$G_1=f_1G_1^*f_1^{-1}$.\endproof

By construction the domain $R(G_1)=R(\Gam'_1)\cap\tau_\pi \left(R(
\Gam'_1)\right)$  is fundamental for the action of $G_1$  in
$\Ome_{G_1}$.

\subheading{Claim 5}  $R(G_1)\cap\Sig =\left(\calP(H_3)\cup
I_{g_1}\cup I'_{g_1}\right)\cap\Sig$.

\pr Recall that $\pi^+(\pi^-)$ means the right (left) component of
${\bf S}^3\bks \pi$\break $\left(I_{g_1}\in \pi^+\right)$.  Then
$\pi^+\cap \Sig\cap R(\Gam'_1)=\calP(H_1)\cap \Sig
=\left(I_{g_1}\cup I'_{g_1}\right)\cap\Sig$  by (b) and (c) of
Lemma 2.

Also, $\tau_\pi\left(\pi^-\cap\Sig\cap\tau_\pi(R(\Gam'_1))
\right)=\pi^+\cap\tau_\pi(\Sig) \cap
R(\Gam'_1)\subset\calP(H_1)\cap\Sig'$, so $\pi^-\cap\Sig \cap
R(G_1)=\tau_\pi\left(\calP(H_1)\right)\cap\Sig =\calP(H_3)\cap
\Sig$.\endprf

Let us consider now the family of spheres $\Sig_t$  centered on
the $y$--axis (figure 2) such that $\Sig_t\cap w_3=\Sig \cap w_3,\
\sig_1=\Sig,\ \sig_0 =\Sig_0,\ t\in [0,1]$,  where $\Sig_t\cap
{\rm ext}(B_1)\cap {\rm ext}(B_2)\subset {\rm ext} (\Sig)\cap {\rm
ext}(B_1)\cap {\rm ext}(B_2)$ (recall ${\rm ext}(\cdot)$ is the exterior
of a set in $\ORR^{\,3}$), $\ \Sig_t\cap I_{g_1}=\emptyset\
(t>0)$. Denote by $\tau_{\Sig_t}$ the corresponding reflections.
As before take the domain $\Ome^*=\Ome_{G_1}\bks G_1 (\Sig^-_0)$
and the group $G'_1={\rm Stab}(\Ome^*,G_1)$, where
$\Sig^-_0={\rm ext}\ (\Sig_0)$  is the unbounded component of
$\ORR^{\,3}\bks \Sig_0$.

Denote $\calG_t=\langle G'_1,\
\tau_{\Sig_t}G'_1\tau_{\Sig_t}\rangle$. Evidently,
$\calG_1=\lim\limits_{t\to 1}\calG_t$.

\proclaim{Lemma 6}  The groups $\calG_t$  are discontinuous, $t\in
[0,1]$.\endproc

\pr First, let us prove the lemma for $t\not= 1$.  By
Claim 5 we have now that $R(G_1)\cap\Sig_t=\calP (H_3)\cap\Sig_t$.
Moreover we claim also that
$$\eqalign{g\Sig_t\cap\Sig_t&=\emptyset ,\ g\in G_1\bks H_3,\ H_3
\Sig_t=\Sig_t, \cr {\rm where}\ \ \ H_3&=\tau_\pi H_1\tau_\pi\
.\cr}\eqno(3)$$

To prove (3) we only need to show that $g(\Sig_t\cap\Lambda
(H_3))\cap (\Sig_t\cap\Lambda (H_3))=\emptyset$,  but this can be
shown from the fact that each point of $\Lambda (H_3)$  is a point
of approximation (see [\P, Claim 1]).

All conditions of Maskit's First  Combination theorem are now
satisfied for the groups $G'_1$ and
$\tau_{\Sig_t}G'_1\tau_{\Sig_t}$ $(t\not= 1)$  [\Mab] and we
obtain also $$\calG_t\cong
G'_1*_{H_3}\left(\tau_{\Sig_t}G'_1\tau_{\Sig_t} \right)\eqno(4)$$
where the $\calG_t$  are all discontinuous, $t\in [0,1)$.

Let us now consider  the group $\calG_1$  and the domain
$R(\calG_1)=R(G_1)\cap \tau_\Sig (R(G_1))$.  Our goal now is to
show that $R(\calG_1)$  is a fundamental domain for the action of
$\calG_1$  in $\Ome_{\calG_1}\left(\infty \in
\Ome_{\calG_1}\right)$.  If now $\langle g_1,\gam_1\nek
\gam_\ell\rangle$  is a set of generators of $G'_1$ then
$S=\langle g_1,\gam_1\nek \gam_\ell,\ g_2,\gam'_1\nek
\gam'_\ell\rangle$  are generators of $\calG_1$, where
$\gam'_i=\tau_\Sig \cdot \gam_i \cdot \tau_\Sig$
 and
$g_2=\tau_\Sig \cdot g_1\cdot \tau_\Sig$. Observe that
 the element $g_1$  is included in $S$ because some
of its isometric spheres belong  to the boundary $\part
R(G'_1)$

We want to apply the Poincar{\'e} Polyhedron theorem [\Mac].
Indeed, an arbitrary cycle of edges in $\part R(\calG_1)$ consists
either of edges situated in $\part (R(G_1))\cap {\rm int} (\Sig)$,
and $\part (\tau_\Sig (R(G_1)))\cap {\rm ext}(\Sig)$,\ \   or
 is an edge
cycle $\ell_1=I_{g_1}\cap I_{g_2},\ \ell_2=I'_{g_1}\cap
I'_{g_2}$, where $I_{g_k},I'_{g_k}$  are the isometric spheres of
$g_k$  and $g^{-1}_k (k=1,2)$.  The sum of angles in any cycle
of the
first type is $2\pi$  because $R(G_1)$  is a fundamental
domain [\Mac].

We now claim that the element $g=g^{-1}_2\cdot g_1$  is parabolic
with a fixed point $d=I_{g_1}\cap I_{g_2}$.  Indeed,
$g_2^{-1}\cdot
g_1=\left(\tau_\Sig\cdot\tau_{I_{g_1}}\right)^2$ because
$g_1=\tau_{\pi_2}\cdot \tau_{I_{g_1}}$  and $\pi_2$  is
orthogonal to $\Sig$ (figure 2).  Now it is easy to check that
$g(d)=d$,  $gI_{g_1}\subset {\rm int} (I_{g_2})$  and
$g({\rm int}(I_{g_1}))={\rm ext}\ (g(I_{g_1}))$,  therefore
the elements
$g$  and $g'=g_1\cdot g\cdot g_1^{-1}$  are parabolics.

All conditions of the Maskit--Poincar{\'e} theorem are valid at the
edges $\ell_i$  also and, hence, $\calG_1$  is discontinuous.
Lemma 6 is proved. \endproof

\proclaim{Lemma 7} The group $\calG_0$  is isomorphic to a
subgroup $L'\subset R$  of a finite index.\endproc

\pr We repeat our construction of $\calG_0$  by modelling it
in ${\bf H}^3$  so as to get the required isomorphism.

Recall that we started from the group $\Gam'_1\subset {\rm
Isom} ({\bf H}^3)$  and showed that
$G_1=\langle\Gam'_1,\tau_\pi\Gam'_1
\tau_\pi\rangle\cong G^*_1=\langle
\Gam'_1,\tau_{w_2}\Gam'_1\tau_{w_2}\rangle$
(see Lemma 4).  Next we constructed $\calG_0$ by using
reflection in $\sig_0=\Sig_0$  such that $\sig_0\cap w_3=\Lambda
(H_3),\ \sig_0\cap B_1= \emptyset,\ w_3=\tau_\pi (w_1)$.

Let $\eta =\tau_{w_2}(w_1)\subset{\bf H}^3,\ \eta\in W$.  Again
let us take the subgroup $G_1^{**}$ of $G_1^*$  which is
$G_1^{**}={\rm Stab}({\bf H}^3\bks G_1^*(\eta^-),\ G^*_1)$, where
$\eta^-$  is a subspace ${\bf H}^3\bks\eta$  not containing $w_2$.

By construction the fundamental domain
$R(G_1^*)=R(\Gam'_1)\cap\tau_{w_2}(R(\Gam'_1))$  of the group
$G^*_1$  satisfies $R(G^*_1)\cap\eta =\calP (H'_3={\rm Stab}
(\eta ,G^*_1))$.  Again by Maskit's First  Combination theorem we have a
group
$L'$:
$$L'=G_1^{**}*_{H'_3}(\tau_\eta G^{**}_1\tau_\eta)\eqno(5)$$
We constructed an isomorphism $\varphi_1\co G^*_1\to G_1$ in Lemma
4 such that $\tau_\pi\cdot\varphi_1\cdot\tau_{w_2}=\varphi_1$,
therefore
$\varphi_1(H'_3)=H_3$  and $\varphi_1(G_1^{**})=G'_1$  .
  It follows now from (4) and (5) that the map
$\varphi_1\big|_{G_1^{**}}$  can be extended to an isomorphism
$\varphi\co L'\to\calG_0$.

Index $|R:L'|$  is finite because $L'$  has a compact
fundamental domain.  The Lemma is proved.  \endproof

Recall that we identify $[\rho]\in \Def(L')$  with
$\rho (L')$.

\proclaim{Lemma 8} There exists a path $\bet_t\co [0,1]\to cl(\Def
(L'))$  such that $\bet_0=L'$, $\bet_1=\calG_1\in \part \Def
(L')$, $\bet_t([0,1))\subset \Def(L')$.\endproc

\pr We have constructed a path $\alpha_t\co [0,1]\to\Def (G^*_1)$ in
Corollary 4 such that $\alpha_0=G^*_1$, $\alpha_1=G_1$ and
$\alpha_t$ is a family of admissible representations.  Let further
$\alpha_t\big|_{G_1^{**}}=\alpha'_t$.  Obviously, the
representations $\alpha'_t$ are also admissible and
$\alpha'_1(G_1^{**})=G'_1$.  We can easily extend our family
$\alpha'_t$  to a family of admissible representations
$\tet_t\co L'\to \Def (L')$ by the formula
$\tet_t=\tau_{\mu_t}\alpha'_t\tau_{\mu_t}$,  where $\mu_t$  are
the spheres constructed in Corollary 4.

Observe that $\mu_1=\pi$  and now take a new continuous family of
spheres $\nu_t$  for which $\nu_t\cap w_3=\Lambda (H_s)=w_3\cap
B_2$  and $\nu_1={\tilde w_3},\ \nu_2=\Sig_0$ where ${\tilde w_3}$
is the sphere containing $w_3$ ($t\in [0,1]$).

Again we have a path $\tet'_t(L')=\langle G'_1,\tau_{\nu_t}
G'_1\tau_{\nu_t}\rangle$.  Composing the path $\tet_t$  with
$\tet'_t$  and with the path corresponding to spheres $\Sig_t$
connecting $\Sigma_0$  with $\Sigma_1$  we get required path
$\bet_t$. The Lemma is proved.
\endproof

\section{Proof of Theorem 1}

(1)\qua Denote by $F=\pi_1\sigma$ a fixed fiber group of our initial
manifold $M$, and let also $F_0=\Gam_0\cap F$.

By J\o rgensen's theorem [\J] the limit $\bet_1=\lim\limits_{t\to
1}\bet_t$  is an isomorphism $\bet_1\co L'\to \calG_1$.   Let us
consider the subgroup $L=L'\cap\Gam_0,\ |\Gam_0:L|<\infty$. Put
also $F_L=L\cap F_0$ for its normal subgroup. We have also the
curve $\bet_t(L)\subset \Def (L)$. Let ${\cal N}=\bet_1(L),\
{\calF}=\bet_1(F_L)$.  Let us show that $g=g_2^{-1}\cdot g_1\in
\calF$. To  this end let us recall that the element $g_1$ was
chosen from the very beginning being in $[H_{w_1}, H_{w_1}]$
(Lemma 1). Recalling also that $\bet_1^{-1}(g_1)=g_1$ and denoting
$\bet_1^{-1}(g_2)=g'_2$,  by construction we get
$g'_2=\tau_{\eta}\cdot g_1 \cdot\tau_{\eta},\ \eta=\tau_{w_2}(w_1)
,\ g_1\in [H_{w_1},H_{w_1}]\subset [F_0,F_0]$  (see Lemma 1).  The
group $\Gam_0$  was chosen to be normal in the reflection group
$R$, and since $[\Gam_0, \Gam_0]\subset F$, it is
straightforward to see that   $$r[F_0,F_0]r^{-1}\subset F_0,\quad
r \in R\ .$$ Hence, $g'_2\in F_0$,  and for the element
$g'=(g'_2)^{-1}\cdot g_1$  we immediately obtain $g'\in
F_L=F_0\cap L'$.  It follows that $\bet_1(g')=g=g_2^{-1}\cdot
g_1\in F_0\cap\calG_1=\calF$ as was promised.

We have that $\cal N$ is isomorphic to the semi-direct product of
$\calF$ and the infinite cyclic group ${\bf Z}$, so taking the
element $t\in \cal N\bks \cal F$ projecting to the generator of
$\cal N/{\calF}$, we observe that the elements $$g_n=t^ngt^{-n}\in
\calF\ ,\quad g\in\calF,\quad n\in {\bf Z}\eqno(6)$$ are all
parabolics. Since $\cal N$ contains no abelian subgroups of rank
bigger than $1$ and $t^n\not\in{\cal F}\ (n\in {\bf Z})$ one can
easily see that the elements (6) are also non-conjugate in
$\calF$. We have proved (1) of the Theorem.

(2)\qua By the construction,  the fundamental polyhedron $R(\calG_1)$
of the group $\calG_1$ contains only one conjugacy class of
parabolic elements $g$  of rank 1. There is a strongly invariant
cusp neighborhood $B_g\cong [0,1]\times R^1\times [0,\infty)$
which comes from the construction of $R(\calG_1)$. So each
parabolic $g_n$ of type (6) gives rise to submanifold
$$B_{g_n}\big/ \langle g_n\rangle\cong T_n\times [0,\infty),\
T_n\cong S^1\times S^1\eqno(7)$$ in the manifold
$M(\calF)=\Ome_N\big/ \calF$. Therefore $M(\calF)$ contains
infinitely many parabolic ends (7) bounded by tori $T_n$. They all
are non-parallel in $M(\calF)$ and therefore by Scott's ``core"
theorem the group $\pi_1(M(\calF))$ is not finitely generated
[\Sc].
\endproof

{\bf Remark}\qua By using the argument of [\P] one can prove:

\proclaim {Theorem 2}  There is a (non-faithful) represention 
$\beta_{1+{\eps}}$ which is $\eps$--close to $\beta_1$ for some
small $\eps>0$ such that the group $\beta_{1+\eps}(F_L)$ is
infinitely generated, has infinitely many non-conjugate elliptic
elements. Moreover, $\beta_{1+\eps}(F_L)$ is a normal infinitely
presented subgroup of a geometrically finite group
$\beta_{1+\eps}(L)$ without parabolics.
\endproc

To prove the theorem one can continue to deform the group for $1<t\leq
1+\eps$ (these representations will no longer be faithful) in order to
get an elliptic element $g_t$ whose isometric spheres form an angle
$\tet(t)$ instead of being tangent. To do this in our Lemma 2, instead
of the sphere $\Sigma$ tangent to the isometric spheres of $g_1$, one needs
to consider a nearby sphere $\Sigma_{1+{\eps}}$ forming angle
$\tet({\eps})$ with them. If $\tet(\eps)={\pi\over 2n}$ and $n>0$ is
large enough the group $\beta_{1+\eps}(F_L)$ is Kleinian, has
infinitely many non-conjugate elliptic elements of the order $n$
(obtained as above as an orbit of $g_{1+\eps}$ by a pseudo-Anosov
automorphism of the $\beta_{1+\eps}(F_L)$). The construction gives us
that $\beta_{1+\eps}(F_L)$ is a normal and finitely generated but
infinitely presented subgroup of the geometrically finite group
$\beta_{1+\eps}(L)$ without parabolic elements. In particular
$\beta_{1+\eps}(L)$ is a Gromov hyperbolic group (see [\P, Lemmas
5--7]).

\references

\Addresses\recd

 \bye